\documentclass[graybox]{svmult}

\usepackage{type1cm}        
\usepackage{makeidx}         
\usepackage{graphicx}        
\usepackage{multicol}        
\usepackage[bottom]{footmisc}

\usepackage{newtxtext}       %
\usepackage{newtxmath}       

\makeindex             

\newcommand{\R}{\mathbb{R}}

\renewcommand{\P}{\mathcal{P}}

\begin{document}

\title*{Recent developments in controlled crowd dynamics}
\author{M. K. Banda \and M. Herty \and T. Trimborn}
\institute{M. K. Banda \at Department of Mathematics and Applied Mathematics, Hatfield 0028, South Africa
 \email{mapundi.banda@up.ac.za}
\and M. Herty \at IGPM, RWTH Aachen, Templergraben 55, 52056 Aachen, Germany  \email{herty@igpm.rwth-aachen.de}
\and T. Trimborn \at IGPM, RWTH Aachen, Templergraben 55, 52056 Aachen, Germany  \email{trimborn@igpm.rwth-aachen.de}
}
%
%
\maketitle

\abstract{We  survey recent results on controlled particle systems. The control aspect introduces new challenges in the 
discussion of properties and suitable mean field limits. Some of the aspects are highlighted in a detailed discussion of a
 particular controlled particle dynamics. The  applied techniques are shown on this simple problem to illustrate the basic
 methods. Computational results confirming the theoretical findings are presented and further particle models are discussed.
 }

\setcounter{tocdepth}{1}

\tableofcontents

\section{Introduction}\label{sec:1}

Large--scale interacting particle systems have recently generated interest in the description of
phenomena beyond classical statistical mechanics and we refer to the books
 \cite{Bellomo:2017aa,Cristiani:2014aa,Pareschi:2013aa} for some examples and further references. A major
 difference to mathematical descriptions in continuum and statistical mechanics is the fact that
 particles are no longer passively interacting. This enables analysis of new pattern formation 
 mechanisms but also allows to introduce control actions within the interacting particle system. 
Among the examples of controlled particle systems are economic models for price formation, wealth accumulation,
 trading or formation of consensus behaviour. We refer to the references \cite{Ben-Naim:2005aa,Aletti:2007aa,During:2009aa,Boudin:2009aa,Wolfram:2011aa,Pareschi:2013aa,Burger:2013ab,Burger:2014ab,Albi:2015aa} as well as
 the references therein for some examples. Typically, control actions might be applied to drive systems towards a desired state
 using either an open loop \cite{Albi:2017aa,Bongini:2015ab,Kalise:2017aa,Herty:2019aa}, a closed loop \cite{Albi:2014ae,Albi:2014ab,Caponigro:2013ab,Degond:2014ab,Albi:2015aa,Albi:2015aa,Herty:2015ab,Albi:2016aa,Albi:2016ab} or a competitive game setup \cite{Huang:2006aa,Huang:2007aa,Huang:2010aa,Lasry:2007aa,Carmona:2013ab,Bensoussan:2013ab,Cardaliaguet:2013ab,Gomes:2014aa,Burger:2014ab,Degond:2014aa,Degond:2014ac,Borzi:2015aa,Cardaliaguet:2017aa,Degond:2017aa,Briceno-Arias:2018aa}. Those directions have been explored recently for interacting particle systems described by ordinary differential equations with  the question of control actions prevailing in the mean field limit of infinitely many particles. Those limits
 have been well--explored in the context of classical statistical mechanics, see, for example, \cite{Cercignani:1988aa,Spohn:1991aa,Cercignani:1994aa}. The limit allows derivations of time--continuous descriptions of qualitative properties independent of the number of particles.  In the case of finitely many interacting agents subject to control, the control problem can be solved using Pontryagin's maximum principle, dynamic programming or Hamilton--Jacobi Bellmann equations. Those techniques are well--established in the context of ordinary differential equations  but pose, in view of the desired mean field limit, interesting and challenging analytical and numerical problems. 
\par The major obstacle, for example, has been the fact that the size of the system of associated Hamilton-Jacobi Bellmann equations is proportional to the number of particles. So far, different approaches have been explored to address possible issues mostly based on an approximation of the control, for example,  by a restriction to short (instantaneous) time horizon controls \cite{Albi:2015aa,Caponigro:2015aa,Herty:2017aa}, or by implementing the control on a binary interaction scale \cite{Albi:2014ab} or by solving approximate Hamilton--Jacobi equations \cite{Albi:2016aa}. Those approaches typically aim to derive a (sub-optimal) closed loop control.  On the level  of  open loop control problems only few and mostly theoretical results could be achieved, for example, in the game theoretic setting \cite{Lasry:2007aa}, using Riccati control \cite{Herty:2015ab} or using a mean field approach for both adjoint and primal variables
   \cite{Herty:2019aa,Albi:2017aa,Fornasier:2014ab}. 
\par Another aspect in dealing with the open loop control problem is a suitable notion of differentiability \cite{Cardaliaguet:2010aa,Bensoussan:2017aa,Herty:2019aa,Herty:2019ab} that leads to consistent optimality systems for both the large--scale particle system and the mean field optimal control problem. Here, choices regarding differentiability in the sense
   of Frechet  as well as differentiability in metric spaces with respect to Wasserstein distance
are possible and an analysis of properties and possible connections is given in \cite[Section 2]{Bensoussan:2017aa} 
or \cite[Appendix]{Herty:2019ab}, amongst many examples. 
\par Suitable numerical methods for the full control problem on the mean field level 
   are still sparse due to the complexity of solving successively forward and backward partial differential equations.  
  Besides the mean field limit, the hydrodynamic limit of  particle systems is of interest for a qualitative study of 
  long--term behaviour of solutions. For controlled particle systems, this direction is currently largely unexplored. 
\par
In this chapter, we would like to illustrate the main steps undertaken in controlled particle systems by presenting, on a simplistic particle
model, challenges and approaches in mean field limits as well as on the hydrodynamic limit. The focus will be on the control 
action and, therefore, we present results in a simple linear setting in Section \ref{sec:2}. This setup has also been used to illustrate basic properties
of mean field games \cite{Cardaliaguet:2010aa} and served as a guiding example for many previously presented techniques \cite{Herty:2019ab}. The linear setting enables use of a closed loop strategy based on the Riccati equation on the level of finitely many agents. 
It is also well--known that in the case of linear--quadratic open loop control problems, the Riccatti control is optimal. Besides
the investigation of the interplay between the Riccati control and the mean field limit, we also discuss here the corresponding hydrodynamic limit in Section \ref{sec:2.3}. The analytical findings are illustrated by numerical simulations of the closed loop control system in Section \ref{sec:3}. Extensions to 
a nonlinear setting and other perspectives are discussed in Section \ref{sec:4}.

\section{Multiscale Riccati Control For Linear Particle Systems}\label{sec:2}

The main purpose of this section is to illustrate control concepts and mean field limits. We,  therefore, restrict
ourselves to a simplistic setting for interacting particle systems. Consider $i=1,\dots,N$ particles
with a state  $(x_i(t),v_i(t)) \in \R^2$  driven by the dynamics 
\begin{align} \label{dyn}
\frac{d}{dt} x_i = v_i, \; \frac{d}{dt} v_i = q_i^*, \; (x_i,v_i)(0)=(x_{i,0},v_{i,0}) 
\end{align}
for initial data $(x_{i,0},v_{i,0})$ and where each particle is subject to a control $q_i^*=q_i^*(t)$ 
modeling  an individual strategy. The control is chosen in order to minimize a joint objective 
\begin{align}\label{ui}
{\boldsymbol q}^*:=\arg\min\limits_{ q \in \R^N} \int_0^T \frac{1}{N} \sum\limits_{i=1}^N \left( \frac{1}2 v_i^2  + \frac{\alpha}2 q_i^2 \right) dt.
\end{align}
where  ${\boldsymbol q}=(q_i)_{i=1}^N$. The parameter $\alpha>0$ is a  weight to balance the cost of all controls and the desired state. The latter is chosen here to 
be $v_i\equiv 0$ for all particles $i=1,\dots,N$ for simplicity. Clearly, other desired states may  be considered and they may 
depedent on  further parameters and costs at terminal 
time $t=T$ may be added to problem \eqref{ui}. The terminal time $T>0$ is fixed but can be replaced by $T=+\infty$ provided  a suitable weight $\exp(- r t )$ is added towards the cost with discount $r>0.$ 

At this point it is important to note that there is a major difference in discussing the case of whether or not all particles try to minimize a joint objective or their individual objective. In the latter
case the problem turns into a differential game and concepts of optimality have to be discussed. We refer to \cite{Friedman:1974aa} for a careful discussion as well as to \cite{Bensoussan:2013ab}, for example, for mean field limits of games. Here, we will focus on the case of  a single joint objective \eqref{ui}.

\subsection{Closed Loop Control For Particle System}\label{sec:particle}
The problem \eqref{ui} and \eqref{dyn} is a linear--quadratic optimal (convex) control problem. It is well--known that in this case the unique  solution ${\boldsymbol q}^*$ is given by a state feedback involving a symmetric matrix $K \in \R^{2N \times 2N}$ which is a solution of the (matrix) Riccati equation. In order to derive this equation, it is advantageous to introduce the following vectors and matrices:
\begin{align} 
&w:=\left( (x_i)_{i=1}^N, (v_i)_{i=1}^N \right)^T \in \R^{2N}, \\
&A: = \begin{pmatrix} 0 & \text{Id} \\ 0 & 0 \end{pmatrix} \in \R^{2N \times 2N}, \; 
B: = \begin{pmatrix} 0 \\ \text{Id}\end{pmatrix}  \in \R^{2N \times N} , \;
M:=  \begin{pmatrix} 0 & 0 \\ 0 &  \text{Id} \end{pmatrix}  \in \R^{2N \times 2N} .
\end{align}
where $\text{Id} \in \R^{N\times N}$ denotes the identity matrix. Solving problem \eqref{dyn}, \eqref{ui} is equivalent to 
solving
\begin{align}
\min\limits_{ q \in \R^N} \int_0^T \frac1{2N} w^T(t) M w(t) + \frac{\alpha}{2N} q(t)^T (\text{Id}) q(t) \; dt \\
\mbox{ subject to }  \frac{d}{dt} w(t) = Aw(t) + Bq(t), \; w(0)=\left( (x_{i,0})_{i=1}^N, (v_{i,0})_{i=1}^N \right)^T 
\end{align}
which necessary optimality conditions yield the following expression for $q\in \R^N$ 
\begin{equation} \label{statefeedback} q(t) = -\frac{2N}{\alpha} B^T K(t) w(t),  \end{equation}
where   $K(t) = K(t)^T$ fulfills the matrix Riccati equation \cite{Sontag:1998aa}:
\begin{align}\label{ricK}
-\frac{d}{dt} K(t) = \frac1{2N} M + K(t) A + A^T K(t) - \frac{2N}{\alpha} K(t) BB^T K(t), \; 
K(T) = 0.
\end{align}
Note that in the nonlinear case a state feedback or closed loop control of the type \eqref{statefeedback} is in general
impossible to obtain and gives rise to  the suboptimal control strategies mentioned in the introduction.
\par 
In view of the mean field limit for $N\to\infty$ particles, it is advantageous to exploit the particular structure of the matrices $A$ 
and $M$ to obtain further information on the matrix $K.$ We split $K$ according to the state space as follows: 
$$K=\begin{pmatrix} K_{11} & K_{12} \\  K_{21} & K_{22} \end{pmatrix},$$
with $K_{ij} \in \R^{N\times N}$ for $i,j=1,2.$ Direct computation yields a simplified form of equation \eqref{ricK} as follows:
\begin{align} \label{ricK2}
-\frac{d}{dt} K(t) =\frac1{2N}\begin{pmatrix} 0 & 0 \\ 0 &  \text{Id} \end{pmatrix} + \begin{pmatrix} 0 & K_{11} \\ 0 & K_{21} \end{pmatrix} + 
\begin{pmatrix} 0 & 0 \\ K_{11} & K_{12} \end{pmatrix} 
- \frac{2N}\alpha \begin{pmatrix} K_{12}K_{21} & K_{12}K_{22} \\ K_{22}K_{21} & K_{22} K_{22} \end{pmatrix}, 
\end{align}
or
\begin{align}
-\frac{d}{dt} {K}_{11} &= - \frac{2N}{\alpha} K_{12}K_{21}, \quad K_{11} (T) = 0,\\
-\frac{d}{dt} {K}_{12} &=K_{11} - \frac{2N}{\alpha} K_{12}K_{22}, \quad K_{12} (T) = 0,\\
-\frac{d}{dt} {K}_{21} &= K_{11} - \frac{2N}{\alpha} K_{22}K_{21}, \quad K_{21} (T) = 0,\\
-\frac{d}{dt} {K}_{22} &= \frac1{2N} Id + K_{12}+K_{21} - \frac{2N}{\alpha} K_{22} K_{22}, \quad K_{22} (T) = 0.
\end{align}
This nonlinear system of ordinary differential equations has a differentiable and locally Lipschitz right hand side. Therefore, 
a unique solution exists and by direct computation, we obtain the assertion of the following Lemma.
\begin{lemma} \label{lem0} 
The system of ordinary differential equations \eqref{ricK} has a unique differentiable solution $K=K(t) \in C^1(\R^{2N \times 2N})$ given by 
  $$K(t)=\begin{pmatrix} 0 & 0 \\  0 & K_{22}(t) \end{pmatrix},$$
  where $K_{22}(t)$ is given as a solution to the equation 
\begin{align}
-\frac{d}{dt} K_{22} (t)= \frac1{2N} \text{Id}  - \frac{2N}{\alpha} K_{22}(t) K_{22}(t), \quad K_{22} (T) = 0.\label{Riccati}
\end{align}
\end{lemma}
Again, the particular structure of equation \eqref{Riccati} allows the derivation of the following result for the structure
of $K_{22},$ that is in fact, diagonal with the same entry. 

\begin{lemma} \label{lem1} The system of ordinary differential equations \eqref{Riccati} has a unique differentiable solution 
\begin{align} \label{sol}
\left( K_{22}(t) \right)_{i,i} = d(t), \mbox{ and } \left( K_{22} (t)\right)_{i,j} = 0, \; i, j=1,\dots,N, \; i \not = j,
\end{align}
where $d(\cdot)$ is the unique solution to the (scalar) ordinary differential equation
$$-\frac{d}{dt} d(t) = \frac1{2N} - \frac{2N}{\alpha} d(t)^2, \quad d(T) = 0.$$
\end{lemma}
The proof of the Lemma follows the verification that $K_{22}$ given by equation \eqref{sol}
 is in fact a solution to equation \eqref{Riccati} which is also the unique solution. The ordinary
 differential equation for $d$ allows for an explicit solution given by 
 $$ d(t) = \frac{1}{2N} \sqrt{\alpha} \tanh\Big( \frac{ T - t}{\sqrt{\alpha} } \Big).$$ However, it is not necessary
 to use the explicit form in the following. Having $d(t)$ available, we obtain  the (optimal) closed loop control 
 \begin{equation}
 {\bf q}^*(t) = - \frac{2N}\alpha B^T K(t) w(t) = - \frac{2N}\alpha (0, K_{22}(t)) w(t) = - \frac{2N}\alpha d(t) v(t). 
 \end{equation}
 Hence, we introduce $y(t):=N d(t)$ that satisfies
 \begin{equation}\label{eqy}
 -\frac{d}{dt} y(t) = \frac1{2} - \frac{2}{\alpha} y(t)^2, \quad y(T) = 0
 \end{equation}
and  $ {\bf q}^*(t) = - \frac{2}{\alpha} y(t) v(t).$ Summarizing, the problem \eqref{dyn}, \eqref{ui} is solved by ${\bf q^*}$ which is the closed loop control. The controlled particle dynamics for $i=1,\dots,N$ is then given by  
\begin{align}\label{ctrl pa}
\frac{d}{dt} x_i = v_i, \; \frac{d}{dt} v_i = - \frac{2}{\alpha} y(t) v_i(t), \; -\frac{d}{dt} y(t) = \frac1{2} - \frac{2}{\alpha} y(t)^2, \\
(x_i,v_i)(0)=(x_{i,0},v_{i,0}), y(T)=0.
\end{align}
The system \eqref{ctrl pa} allows for a mean field limit in the number of agents as shown in the following paragraph.
 
 Furthermore, we note that the decay of $v_i(t)$ towards the desired zero state can be quantified using a Lyapunov function. The following result
 holds true. 
 \begin{lemma} \label{lemP} 
 Consider $N$ particles with dynamics given by equation \eqref{ctrl pa} and arbitrary initial conditions $(x_{i,0},v_{i,0})$ for $i=1,\dots,N.$ 
 Let $y$ be the solution to \eqref{eqy}. 
 Then, the differentiable function 
 $$L(t):= w^T(t) K(t)w(t)$$
 where $K$ is given by equation \eqref{ricK}  is bounded from above as 
 $$L(t) \leq L(0) \exp\left( -  r (t) \right)$$
 for  $t\in [0,T]$ and the rate $r(t)$ is given by 
 \begin{equation}\label{rate}
 r(t)=   \int_{0}^t \frac{2y(s)}\alpha ds \geq 0.
 \end{equation}
 \end{lemma}
 \noindent {\bf Proof}.  We observe that due to Lemma \ref{lem1} $d(t)\geq 0$ for $t \in [0,T]$ and therefore $K(t)$ is positive semi--definite. Hence,  $L(t)\geq 0$ for all $t\in  [0,T].$ According to Lemma \ref{lem0} and Lemma \ref{lem1}, we have $ L(t) =  \sum\limits_{i=1}^N v_i(t)^2  d(t) = \frac1N \sum\limits_{i=1}^N v_i^2(t) y(t). $ Thus
 \begin{align*}
\frac{d}{dt} L(t) &= \frac{1}N \frac{d}{dt} y(t) \sum\limits_{i=1}^N v_i(t)^2  -  \frac{4}{N \alpha}  y^2(t) \sum\limits_{i=1}^N v_i^2(t) \\ 
&= - \frac12   \frac1N \sum\limits_{i=1}^N v_i(t)^2  + \left( \frac{2}{\alpha } -  \frac{4}{ \alpha} \right) y(t) L(t) \leq -\frac{2 y(t) }{\alpha} L(t). 
 \end{align*}
The  assertion follows by Gronwall's inequality and since $y \geq0$ the rate is non--positive. 
 
 \subsection{Mean Field Limit of Controlled Particle System}

We illustrate the mean field limit of the particle system \eqref{ctrl pa}. For more details and a detailed discussion 
of the mean field limit for linear systems, we  refer to \cite{Golse:2016aa,Jabin:2014aa}. The following formal calculation 
exhibits the main idea: Denote by $\mu^N(t,\cdot) \in \P(\R^2)$ the empirical measure  at time $t$ associated with the 
 state $( x_i(t),  v_i(t))_{i=1}^N \in \R^{2N}$  
by 
\begin{equation}\label{empirical}
\mu^N(t,x,v) := \frac1N \sum\limits_{i=1}^N\delta(x-x_i(t)) \delta(v-v_i(t)),
\end{equation}
where $\delta$ denotes the Dirac measure on $\R$ and $\P(\R^2)$ is the space of probability measures on $\R^2.$   Let $\psi$ be
 any smooth  compactly supported function, i.e.,  $\psi(x,v)\in C^{\infty}_0(\R^2)$. In the following formal computation, we neglect
 the time-dependence of the corresponding functions for readability.
\begin{align*}
\frac{d}{dt} \int_{\R^2} \psi(x,v) d\mu^N(x,v)&= \frac1N \sum\limits_{i=1}^N \partial_x  \psi(x_i,v_i)\frac{d}{dt}{x}_i +\partial_v \psi(x_i,v_i)\frac{d}{dt} {v}_i \\
&= \frac1N \sum\limits_{i=1}^N \partial_x  \psi(x_i,v_i) v_i  - \frac{2}\alpha y \partial_v \psi(x_i,v_i) v_i \\
&= \int_{\R^2} \partial_x \psi(x,v) v d\mu^N(x,v) - \frac{2}\alpha y \int_{\R^2} \partial_v \psi(x,v) v d\mu^N(x,v)
\end{align*}
Hence, the measure $\mu^N$ fulfills a partial differential equation in the sense of distributions. In the following, we will assume 
that  the measure $d\mu(t,x,v)$ has a density $f=f(t,x,v)$, i.e., $d\mu(t,x,v)=f(t,x,v) dxdv,$ then 
the previous equality is the weak form of the mean field equation for $(x,v) \in \R^2$ and $t\in [0,T]$
\begin{equation}\label{mf}
\partial_t f(t,x,v) + \partial_x  \left( v f(t,x,v) \right) - \frac{2}\alpha y(t) \partial_v \left( v f(t,x,v) \right) = 0.
\end{equation}
Equation \eqref{mf} has to be solved subject to  initial conditions 
\begin{equation} f(0,x,v) = f_0(x,v) \end{equation} 
where the initial (non--negative) 
probability density $f_0(x,v)$ is an approximation to the empirical measure $\mu^N_0$ associated with the initial conditions $\mu^N_0=\frac{1}N 
\sum\limits_{i=1}^N\delta(x-x_{i,0}) \delta(v-v_{i,0}).$  Furthermore, $y$ fulfills, as before, 
 \begin{equation}\label{eqy2}
 -\frac{d}{dt} y(t) = \frac1{2} - \frac{2}{\alpha} y(t)^2, \quad y(T) = 0. 
 \end{equation}
\par
In order to show that the particle dynamics converges to the mean field limit, the Wasserstein distance and Dobrushin's inequality has been 
used as a theoretical tool. We follow the presentation in \cite{Golse:2016aa} and references therein for further details. The convergence 
is obtained in the space of probability measures  $\P(\R^2)$ using the Wasserstein distance. This distance measures the space of
probability measures and we refer to \cite{Ambrosio:2008aa} for more details. For the following presentation it sufficies to consider the 
$1-$Wasserstein distance defined as follows.
\begin{definition} Let $\mu$ and $\nu$ be two probability measures on $\R^2.$ Then, the $1-$Wasserstein distance is defined 
by \begin{align} W(\mu,\nu):= \inf\limits_{ \pi \in \P^*(\mu,\nu) } \int_{\R^2}\int_{\R^2} | \xi - \eta | d\pi(\xi,\eta) \end{align}
where $\P^*(\mu,\nu)$ is the space of probability measures on $\R^2 \times \R^2$ such that the marginals of $\pi$ are $\mu$ and $\nu,$ respectively, i.e., 
$\int_{\R^2} d\pi(\cdot,\eta) = d\mu(\cdot)$  and $\int_{\R^2} d\pi(\xi,\cdot) = d\nu(\cdot)$.
\end{definition}
Furthermore, we introduce the push--forward notion for a measureable map $g:\R^2 \to \R^2$ and a measure $\mu \in \P(\R^2)$. A measure $\nu \in \P(\R^2)$ 
is denoted by $ \nu = g \# \mu,$ if 
\begin{align*} \nu(A) = \mu( g^{-1}(A) ) \end{align*}
for any set $A \subset \R^2$ such that $g^{-1}(A)$ is $\mu-$measurable. Let $y$ be the unique solution to equation \eqref{eqy2}.  Introduce now $\xi=(x,v) \in \R^2$ and rewrite  
equation \eqref{mf}  for a probability measure $\mu(t,\xi):=f(t,x,v)dx dv $ as
 \begin{equation}\label{mfxi} \partial_t \mu(t,\xi) + \nabla_\xi \cdot \left( ( \xi_2, - \frac{2}\alpha y(t) \xi_2)^T \mu(t,\xi) \right) =0, \; \mu(0,\xi)=\mu_0:=f_0(x,v)dxdv.
 \end{equation}
 which is understood in the sense of distributions.  Associated with the mean field equation \eqref{mfxi} and the initial datum, there is a system of characteristics  $t \to \Xi(t, \bar{\xi}) \in \R^2$  emanating from
 $\bar{\xi} \in \R^2$  by 
 \begin{align}\label{ode11}
\frac{d}{dt} \Xi_1(t, \bar{\xi}) &= \Xi_2(t, \bar{\xi}), \;  \Xi_1(0, \bar{\xi})=\bar{\xi}_1, \\
\frac{d}{dt} \Xi_2(t, \bar{\xi}) &= - \frac{2}\alpha y(t) \Xi_2(t, \bar{\xi}), \;  \Xi_2(0, \bar{\xi})=\bar{\xi}_2.\label{ode12}
 \end{align}
For any $t\geq 0$ the measure $\mu(t,\cdot)$ is then obtained as 
\begin{align} \label{lagrange} \mu(t,\cdot)= \Xi(t,\cdot) \# \mu_0.\end{align}
The latter equation follows by integral transformation formula applied to equation \eqref{mfxi}. The relation \eqref{lagrange} is
now the starting point for developing the Dobrushin estimate. Let  measures $\mu=\mu(t,\cdot)$ and $\nu=\nu(t,\cdot)$  in $\P(\R^2)$ 
be two solutions to \eqref{mfxi} such that  $t\to\mu(t,\cdot)$ and $t\to \nu(t,\cdot)$ are continuous. Then, we have
\begin{align}\label{d ineq}
W(\mu(t,\cdot), \nu(t,\cdot) )\leq W(\mu_0,\nu_0) C_1, \; t \in [0,T]
\end{align}
for some constant $C_1\geq 0.$ The inequality \eqref{d ineq} is established using the following computation. Let $\pi_0 \in \P^*(\mu,\nu)$ and for $(\xi,\eta)\in\R^4$ define $\pi(t,\cdot)$ as the measure under the image of $(\xi,\eta)\to(\Xi(t,\xi),\Xi(t,\eta))$ by 
$$\pi(t, (\xi,\eta) ):=\pi_0( \Xi(t,\xi), \Xi(t,\eta) ).$$ Then, we have $\pi(t,\cdot) \in \P^*\left( \mu(t,\cdot),\nu(t,\cdot) )\right).$ 
Further, denote by $R(t)$ the distance computed at the measure $\pi(t,\cdot)$ as 
\begin{align*}
R(t):=\int_{\R^4} \|\xi-\eta\| d\pi(t,(\xi,\eta)) = \int_{\R^4}  \|\Xi(t,\xi)-\Xi(t,\eta)\| d\pi_0(\xi,\eta) \\
=\int_{\R^4} \| \xi - \eta + \int_0^t (\Xi_2(s,\xi) -\Xi_2(s,\eta) ,- \frac{2}\alpha y(s) ( \Xi_2(s,\xi) - \Xi_2(s,\eta) )^T ds \| d\pi_0(\xi,\eta).
\end{align*}
Note that the solution to \eqref{eqy2} fulfills $0\leq y(t)\leq \frac{\sqrt{\alpha}}{2}$ for $t \in [0,T].$ Furthermore, the right-hand side
of the system \eqref{ode11}-\eqref{ode12} is Lipschitz with constant $C_L:=\max\{1 , \frac{1}{\sqrt{\alpha}} \}$ and therefore we have,
for any fixed $T$, using Gronwall's inequality:
\begin{align} 
\| \Xi(s,\xi) -\Xi(s,\eta)\|^2 \leq \exp( C_L T ) \| \xi-\eta\|^2, \; \forall s \in [0,T].
\end{align}
Hence, we estimate 
\begin{align*}
R(t)&\leq \int_{\R^4} \left( \| \xi-\eta \| + 4 T \exp( \frac{1}2 C_L T ) \| \xi-\eta\| \right) ds d\pi_0(\xi,\eta) \leq C_1 R(0)
\end{align*}
for a constant $C_1=\max\{ 1, 4 T \exp( \frac{1}2 C_L T ) \}.$ Next, we take the infimum for both sides over all $\pi_0 \in \P^*(\mu_0,\nu_0)$ 
and obtain the Dobrushin type inequality\eqref{d ineq}, i.e., 
\begin{align*}
W(\mu(t,\cdot), \nu(t,\cdot) ) \leq  R(t) \leq C_1 W(\mu_0,\nu_0).
\end{align*}
This inequality is a key estimate to show convergence of the particle dynamics towards the mean field equation: consider
$\nu_0^N$ to be the empirical measure associated with initial data $\xi_{0,i}=(x_{0,i},v_{0,i})$ for $i=1,\dots, N,$ i.e.
\begin{align} \nu_0^N(\xi) := \frac{1}N \sum\limits_{i=1}^N \delta( \xi - \xi_{0,i}). \end{align}
Then, \eqref{lagrange} implies that the measure at time $t$ is given by the empirical measure
\begin{align} \nu^N(t,\xi) := \frac{1}N \sum\limits_{i=1}^N \delta( \xi - \Xi^{i}(t)), \end{align}
where $\Xi^{i}(t)= \Xi(t, \xi_{0,i})$ is the solution to the system \eqref{ode11}--\eqref{ode12} with initial data $\xi_{0,i}=(x_{0,i},v_{0,i}),$ i.e., 
 \begin{align}\label{ode21}
\frac{d}{dt} \Xi^{i}_1(t) &= \Xi^{i}_2(t), \;  \Xi^{i}_1(0)=x_{i,0}, \\
\frac{d}{dt} \Xi^{i}_2(t) &= - \frac{2}\alpha y(t) \Xi^{i}_2(t), \;  \Xi^{i}_2(0) =v_{i,0}.\label{ode22}
 \end{align}
Recalling that $\xi=(x,v),$ the system \eqref{ode21}--\eqref{ode22} is equivalent to the controlled particle system \eqref{ctrl pa}. 
 Let $\mu=\mu(t,\xi)$ be  another solution to the mean field equation \eqref{mfxi}. Then,   Dobrushin's inequality shows
that 
$$W(\mu(t,\cdot), \nu^N(t)) \to 0$$
provided that $W(\mu_0,\nu_0^N)\to 0.$ Hence, we have shown that the limit of the particle system converges in Wasserstein
towards a weak solution to the mean field equation.
\par
Finally, we also prove the decay of the mean field Lyapunov function by extending the results of Lemma \ref{lemP}. For the empirical measure \eqref{empirical}  and its density, $f$,  we obtain as in the proof of Lemma \ref{lemP}.
\begin{align}
L(t) &= w^T K(t) w=   \frac1N \sum\limits_{i=1}^N v_i^2(t) y(t) = \int_{\R^2} v^2 y(t) d\mu^N(t,x,v) \\
&= \int_{\R^2} v^2 y(t) f(t,x,v) dxdv .
\end{align}
The following Lemma shows that the mean field limit inherits  the decay properties  of the particle system and the same rate. 
\begin{lemma} \label{lemF} 
Consider sufficiently smooth initial conditions $f_0(x,v) \in \P(\R^2)$ with 
 integrable second moment  $\int_{\R^2} v^2 f_0(x,v) dxdv  < \infty.$ Assume that there exists 
  a non--negative $C^1([0,T]\times \R^2)$ solution $f=f(t,x,v)$ vanishing as $\|(x,v)\|\to \infty$ 
 to the dynamics given by equation \eqref{mf}  with 
 integrable second moment  $\int_{\R^2} v^2 f(t,x,v) dxdv  < \infty.$ Let $y$ be the solution to \eqref{eqy2}. 
 Then, the differentiable function 
 $${\bf L}(t):= \int_{\R^2} v^2 y(t) f(t,x,v) dxdv $$
  is bounded from above by 
 $${\bf L}(t) \leq {\bf L}(0) \exp\left( -  r(t)  \right)$$
 for  $t\in [0,T]$ and  rate, $r(t)$, given by equation \eqref{rate}.
 \end{lemma}
 \noindent {\bf Proof}. The (strong) assumptions allow for a pointwise evaluation of the partial differential equations. Furthermore, since $y\geq 0$ and $f$ is assumed to be non--negative, ${\bf L}\geq 0.$ A direct computation of the derivative of ${\bf L}$ yields 
 \begin{align*}
 \frac{d}{dt}{\bf L}(t) &= \int_{\R^2} v^2 y(t) \left( - \partial_x ( v f(t,x,v) )  + \frac{2}\alpha y(t) \partial_v ( v f(t,x,v) ) \right) dxdv + \\ 
 & + \frac{d}{dt} y(t) \int_{\R^2} v^2  f(t,x,v) dxdv  \\
 \leq &   -  y(t) \int_{\R}  \partial_x \left( \int_{\R} v^3  f(t,x,v) dv \right) dx - \int_{\R^2} \frac{2}\alpha y^2(t) 2 v^2 f(t,x,v) dx dv    \\
&  + \frac{2}\alpha y(t) {\bf L}(t), 
 \end{align*}
where we used integration by parts in the $v$ derivative and, since $f$ is non--negative,  could estimate $\frac{d}{dt}y$ using equation \eqref{eqy2} .
Further, we obtain 
 \begin{align*}
 \frac{d}{dt}{\bf L}(t) &\leq \left( -\frac{4}{\alpha} + \frac{2}\alpha \right)y(t)  {\bf L}(t) = -\frac{2}\alpha y(t)  {\bf L}(t). 
 \end{align*}
Using Gronwall's inequality  the same rate $r(t)$ as in Lemma \ref{lemP} is obtained.  This shows that the basic properties of the controlled system are transferred from particle to mean field limit. 

 \subsection{Hydrodynamic Approximation to Controlled Particle System}\label{sec:2.3}
 
As the final scale of description, we consider a hydrodynamic formulation of the particle system. The procedure to derive hydrodynamic equations 
from kinetic or mean field equations is by now standard and a detailed discussion of the validity and its properties could be found in the references such as \cite{Cercignani:1994aa,Levermore:1996aa}. In view of equation \eqref{mf}, we define the quantities $\rho(t,x)$ as local particle density as well as the flux $(\rho u)(t,x)$ at time $t$ and 
position $x$ by
\begin{align}\label{macro}
\rho(t,x):=\int_{\R} f(t,x,v) dv, \; 
(\rho u)(t,x):=\int_{\R} v f(t,x,v) dv. 
\end{align}
Here, $u=u(t,x)$ is the hydrodynamic velocity. Further, we introduce the pressure $p[f]$ as 
\begin{align} p[f](t,x) = \int_{\R} (v-u(t,x))^2 f(t,x,v)dv.\end{align}
Formally, evolution equations for the density and flux are found by integration of the mean field equation \eqref{mf} with respect to $dv$ and $v dv.$ 
The quantities \eqref{macro} are also referred to as the zeroth and first moments of the probability density $f$ with respect to $v.$ The partial differential equations
obtained by integration with respect to $x \in \R$ and $t \in[0,T]$ are given by 
\begin{align}\label{hydro}
\partial_t \rho(t,x) + \partial_x (\rho u)(t,x) = 0, \\ 
\partial_t (\rho u)(t,x) + \partial_x \left( (\rho u^2)(t,x) + p[f](t,x) \right) + \frac{2}\alpha y(t) (\rho u)(t,x) = 0.\label{hydro2}
\end{align}
The equation for $y$ is unchanged and given by equation \eqref{eqy2}. Initial conditions for the previous system are obtained by successive integration 
of the initial condition, i.e., 
\begin{align}\label{hydro-IC}
\rho(0,x)= \rho_0(x) := \int _{\R} f_0(x,v)dv, \; (\rho u)(0,x)=(\rho u)_0(x):=\int _{\R} v f_0(x,v)dv.
\end{align}
The system \eqref{hydro}-\eqref{hydro2} is not closed  since $p$ can not be expressed explicitly in terms of $(\rho, \rho u)$ directly. This problem is known as a closure problem and it is due to the fact,
that $f$ is projected to a lower dimensional subspace given by density and flux and $p$ accounts for the projection error. Several possibilites exist to overcome the problem. Further equations for the moments of the type $m^j(t,x)=\int_{\R} v^j f(t,x,v) dv$ of $f$ can be derived shifting the closure problem to equations for higher moments and thereby improving the approximation of $f$ in terms of $(\rho, \rho u, m^3,\dots)$. However, the problem in principle persists. A possibility to close the hierarchy of equations is to assume that $f$ is close to a (known) equilibrium distribution  for higher moments. We illustrate this approach by approximating 
\begin{equation} \label{eqP} p[f] \approx p[ \delta(\cdot - u(t,x) )]\end{equation}
known as mono-kinetic closure. Equation \eqref{eqP} states that at the level of the pressure the particles propagate with the hydrodynamic speed, instead
of their individual velocity $v.$ Formally, replacing $p[f]$ by the mono--kinetic closure in equation \eqref{hydro}--\eqref{hydro2} leads to the system similar to pressureless gas dynamics of the form: 
\begin{align}\label{hydro-pre}
\partial_t \rho(t,x) + \partial_x (\rho u)(t,x) = 0, \\ 
\partial_t (\rho u)(t,x) + \partial_x  (\rho u^2)(t,x)  + \frac{2}\alpha y(t) (\rho u)(t,x) = 0.\label{hydro-pre2}
\end{align}
Another used closure in gas dynamics is the Grad--closure \cite{Cercignani:1994aa}. Here, we
assume that 
\begin{equation} \label{eqG} p[f] \approx C_1 \rho^{C_2} \end{equation}
for some known constants $0\leq C_1,C_2.$ In the case of isentropic gas $1< C_2 \leq 3$ and for isothermal gas, we have $C_2=1.$ Other closure options, for example, based on entropy principles also exist, see e.g. \cite{Dubroca:1999aa}.
\par
Based on the moments of the mean field $f$, we obtain 
\begin{align}
{\bf L}(t) = \int_{\R^2} v^2 y(t) f(t,x,v) dxdv = \\
y(t) \int_{\R} \left( \int_{\R} (v-u(t,x))^2 f(t,x,v) dv \right) +  (\rho u^2)(t,x) dx  = \\ y(t) \int_{\R} p[f](t,x) + (\rho u^2)(t,x) dx.
\end{align}
For the closure \eqref{eqP}, we are able to also prove decay of the corresponding Lyapunov function at the hydrodynamic scale. 
\begin{lemma} \label{lemH} 
Consider sufficiently smooth initial conditions $\rho_0(x), (\rho u)_0(x)$ such that  
 $\int_{\R} (\rho u^2)_0(x) dx  < \infty.$ Assume that there exists a
   $C^1([0,T]\times \R)$ solution $(\rho,\rho u)(t,x)$  vanishing at $\|x\|\to \infty$  and with $\rho \geq 0$
 to the dynamics given by equation \eqref{hydro-pre}-\eqref{hydro-pre2}  with 
 integrable momentum  $\int_{\R} (\rho u^2)(t,x) dx  < \infty$ for any time $t.$ Let $y$ be the solution to \eqref{eqy2}. 
 Then, the differentiable function 
 $${\mathcal L}(t):=  \int_{\R} y(t) (\rho u^2)(t,x)  dx $$
  is bounded from above by 
 $${\mathcal L}(t) \leq {\mathcal L}(0) \exp\left( -  r(t)  \right)$$
 for  $t\in [0,T]$ and  the rate, $r(t)$, is given by equation \eqref{rate}.
 \end{lemma}
 \noindent {\bf Proof}. The (strong) assumptions allow for a pointwise evaluation of the partial differential equations. Note that equations
 \eqref{hydro-pre} and \eqref{hydro-pre2} lead to an equation for $u=u(t,x)$ as 
 \begin{align*} \partial_t u(t,x) + u(t,x) \partial_x u(t,x) + \frac{2}\alpha u(t,x) = 0. \end{align*}
 Furthermore, a direct computation yields
 \begin{align*} 
\frac{d}{dt} {\mathcal L}(t)& =  \frac{d}{dt} y(t)  \int_{\R} (\rho u^2)(t,x)  dx  +  \\
& y(t) \int_{\R} (-  \partial_x (\rho u)(t,x) )  u^2(t,x) dx + 2 (\rho u)(t,x) \partial_t u(t,x) dx  \\
\leq & \frac{2}{\alpha} y^2(t)  \int_{\R} (\rho u^2)(t,x)  dx  + y(t) \int_{\R} 2 (\rho u)(t,x)  \left( u(t,x) \partial_x u(t,x) + \partial_t u(t,x) \right) dx \\
\end{align*}
  where we used the equation for $y$ and the non--negativity of $\rho(t,x)$. Further, using the equation for $u(t,x)$, we obtain 
    \begin{align*} 
\frac{d}{dt} {\mathcal L}(t)& \leq \frac{2}{\alpha} y(t) {\mathcal L}(t) - y(t)  \int_{\R}    2 (\rho u)(t,x) \frac{2}\alpha u(t,x) dx = \frac{-2 y(t)}\alpha {\mathcal L}(t). 
\end{align*}
Using Gronwall's inequality we obtain the assertion. The obtained rate is the same as on the particle and mean field scale. 

 \begin{remark} Using Grad's closure a similar estimate is not possible due to the particular structure of the hydrodynamic equations
 and the form of the Lyapunov function.
 \end{remark}



\section{Numerical Tests and Results for the Linear Control Problem} \label{sec:3}

Numerical experiments will be undertaken to demonstrate the theoretical properties discussed in the previous sections. The first set of numerical tests will consider the particle system. The theoretical results are presented in Section \ref{sec:particle}. The decay of Lyapunov function will be investigated and presented in Section \ref{sec:numparticle} below. Similarly, the hydrodynamic limit discussed in Section \ref{sec:2.3} will be numerically investigated and presented in Section \ref{sec:numhydro}. 

\subsection{Numerical results for the particle system}\label{sec:numparticle}

In these tests, we consider a system with a total of $N=250$ particles. The initial spatial positions of these particles,  $x_{i,0}$, are uniformly distributed in the interval $[0,1].$ The initial velocities are given by $v_{i,0}=\exp(x_{i,0} \sin(2 \pi x_{i,0} ) + \xi_i$ where $\xi_i$ is a random variable uniformly distributed in $[0, 0.2].$ The dynamics are described by equation \eqref{ctrl pa}. Time marching of this dynamic system is undertaken using an embedded explicit Runge--Kutta method of orders 2 and 3 with four stages (with local error control). The control $y(t)$ is computed by piecewise constant discretization of the exact solution to equation \eqref{eqy2} on 
the same grid used for the integration of the system of ordinary differential equations \eqref{ctrl pa}.

The Lyapunov function $L(t)$ in Lemma \ref{lemP} is also evaluated on the same grid.  For comparisons the theoretical decay rate, $r(t)$, is also computed by a midpoint integration of equation \eqref{rate} using the exact solution to the ordinary differential equation \eqref{eqy2}. The terminal time for these computations has been set to $T=1.$

\begin{figure}[htb]
\center \includegraphics[width=\textwidth]{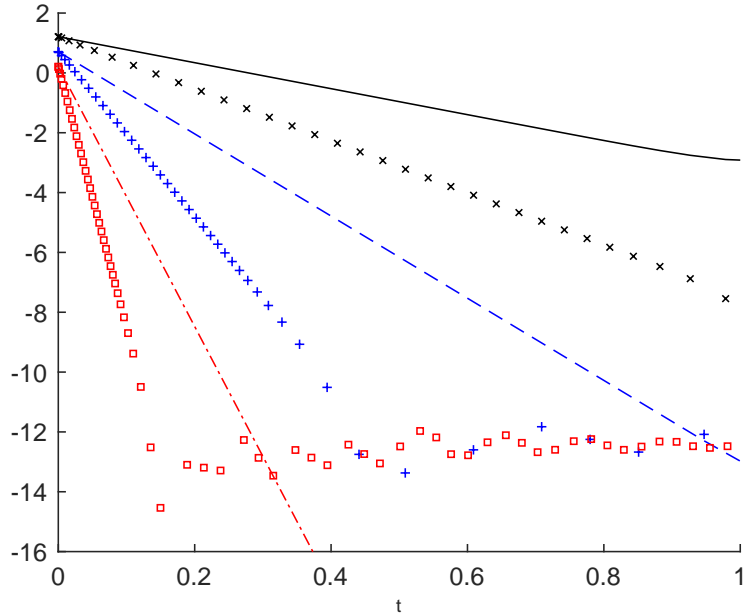}
\caption{Semi--logarithmic plot for the decay of the Lyapunov function $t\to L(t)$ over time. The crossed lines correspond
to the numerical integration of system \eqref{ctrl pa} for different values of $\alpha.$ The solid line shows the theoretical 
expected decay $L(0)\exp(-r(t))$ where the rate $r$ is obtained by numerical integration of equation \eqref{rate}. 
The red,blue and black line correspond to $\alpha=10^{-4},$ $\alpha=10^{-3}$ and $\alpha=10^{-2},$ respectively. }
\label{fig1}
\end{figure}

In Figure \ref{fig1} computational results for the particle system are presented. The decay of the Lyapunov function, $L(t)$, with the estimate of the decay provided by Lemma \ref{lemP} for different values of $\alpha$ is computed. We consider $\alpha \in \{ 10^{-2}, 10^{-3}, 10^{-4} \}.$ Recall that $\alpha$ weighs the control cost compared to the cost associated with the quadratic deviation of particle velocities from zero. Larger values of $\alpha$ are related to higher control costs and lead to slower decay rates which is also observed numerically. Furthermore, we observe that the theoretical estimate is an upper bound on the observed decay in all cases. This confirms the theoretical findings.

\subsection{Numerical results for the hydrodynamic limit}\label{sec:numhydro}
Similar computations as above are repeated for the hydrodynamic approach in Section \ref{sec:2.3}. Note that formally the equations are the so-called pressureless gas dynamics - however damped by the control in $\rho u.$ Without the control the equations may exhibit Dirac-solutions for initial values in the velocity field that lead to concentration phenomena. This case
is not present for the initial data considered here due to the strong damping and the smooth initial datum. The initial data of the particle system translates to the following initial data for the hydrodynamic simulation 
\begin{equation} \rho_0(x) = 1 \mbox{ and }  (\rho u)_0(x) = \exp(x) \sin(2\pi x) + \xi(x) \end{equation}
where $x\to\xi(x)$ is spatially distributed random noise. Periodic boundary conditions in space  are imposed. 

A second--order relaxed finite--volume scheme for the spatial and temporal discretization of equation \eqref{hydro} as proposed in \cite{Jin:1995aa} is employed. The spatial domain is discretized using an equal-distant grid with center points $x_i = i \Delta x$ with $i=1,\dots,N_x$ grid points. We choose  $N_x=250$ spatial points.   The temporal grid is chosen using
a CFL condition with a CFL number of $0.9.$ Terminal time is set to $T=1.$ The noise is pointwise equally distributed in $[0, 0.2].$ As in the case
of the particle system the control $y(t)$ is computed by  piecewise constant discretization of the exact solution to equation \eqref{eqy2} on 
the grid used for the solution of partial differential equations \eqref{hydro}. 

It is appropriate to remark at this point that in \cite{Bouchut:2004aaa} a first-order version of the relaxed finite-volume scheme also referred to as the Rusanov scheme was discussed. For a system of conservation laws:
\[\partial_t U + \partial_x F(U) = 0\] the numerical flux takes the form:
\[F_{i+1/2} = \frac{F(U_i) + F(U_{i+1})}2 - c \frac{U_{i+1} - U_i}2\]
The subcharacteristic condition needs to be satisfied i.e. the eigenvalues of the physical system to be solved must lie between the eigenvalues of the relaxation system. It was proved that for isentropic gas dynamics systems with pressure $p = p(\rho)$ the choice of 
 \[c = \max\Big(|u_i| + \sqrt{p^\prime(\rho_i)}, |u_{i+1}| + \sqrt{p^\prime(\rho_{i+1})}\]   
where $u$ is the velocity gives a scheme that preserves positivity of density as well as resolves solutions with vacuum since $c$ does not blow up at vacuum. For relaxation schemes that can handle delta-shocks, we refer to \cite{Berthon:2006aaa}. For other numerical approaches for pressureless gas, we refer to \cite{Bouchut:2003aaa,Boudin:2012aaa,Chalons:2012aaa}.

The Lyapunov function $\mathcal{F}(t)$ 
in Lemma \ref{lemH} is evaluated on the same temporal grid and the integral in space is discretized using a midpoint scheme.

\begin{figure}[htb]
\center \includegraphics[width=\textwidth]{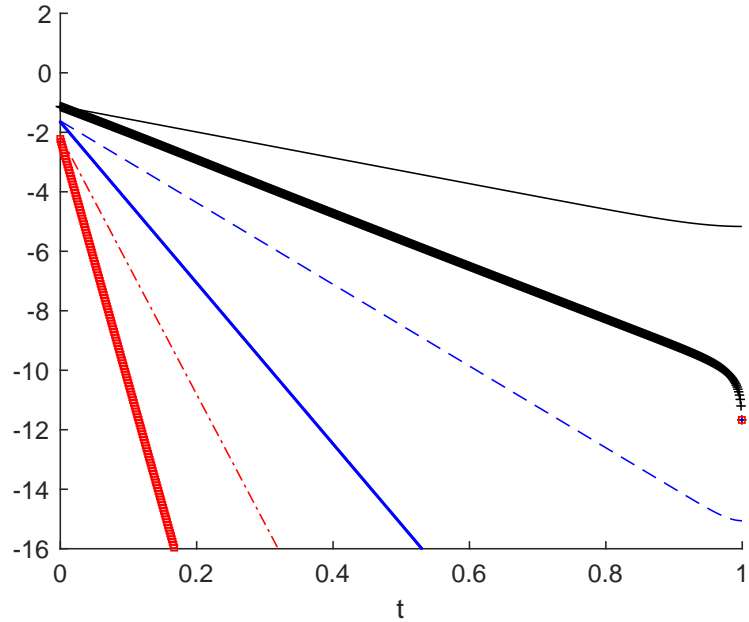}
\caption{Semi--logarithmic plot for the decay of the Lyapunov function $t\to {\mathcal L}(t)$ over time. The crossed lines correspond
to the numerical integration of system \eqref{ctrl pa} for different values of $\alpha.$ The solid line shows the theoretical 
expected decay ${\mathcal L}(0)\exp(-r(t))$ where the rate $r$ is obtained by numerical integration of equation \eqref{rate}. 
The red,blue and black line correspond to $\alpha=10^{-4},$ $\alpha=10^{-3}$ and $\alpha=10^{-2},$ respectively. }
\label{fig2}
\end{figure}

In Figure \ref{fig2}, we compare the decay of the Lyapunov function $\mathcal{L}(t)$ with the estimate of the decay provided by Lemma \ref{lemH}
for different values of $\alpha.$ As in the case of the particle system, the theoretical findings are confirmed by the numerical simulation.

\section{Results on Nonlinear Interacting Particle Systems} \label{sec:4}

In this section we briefly review some existing research directions for controlled particle systems. 

\subsection{One--dimensional State Space Models}
Many recent publications focus on simple models with a phase state that is one--dimensional, i.e., 
each particle $i$ has a state $y_i \in \R$.  Such models are popular in the description of opinion formation and wealth distributions \cite{Hegselmann:2002aa,Naldi:2010aa} to name just two.
\par
Many contributions \cite{Albi:2014ae,Albi:2016ab,Albi:2015aa} introduce novel control strategies based on opinion formation models of the type 
\begin{align}\label{y}
&y_i' = \frac{1}N \sum\limits_{j \not = i }^N P(y_i,y_j) (y_j-y_i) + q_i, \; y_i(0)=y_{i,0}.
\end{align}
where $y_i \in I \subset \R$ and the dynamics are driven by an alignment process due to  pairwise interaction and weighted by a function $P(\cdot,\cdot)$. In the case of wealth models, we have $I=\R_+$ denoting the available
money for each particle while for opinion formation, we have $I=[-1,1]$ denoting extreme opinions. 
For traffic flow applications the interval $I$ is $I=[0,v_{\max}]$ where $v_{\max}$ is the maximal speed allowed 
on a road. Other examples are $I=[0,2 \pi]$ in the case of one-dimensional Kuramoto--type models or $I=\R$ 
in the case of the one--dimensional Cucker--Smale model. 

For constant $P$, the model \eqref{y} is still linear 
and a similar analysis as in Section \ref{sec:2} is possible. For nonlinear $P$, most of the research has focused
on suboptimal control based on instantaneous or short time horizon control. The mean field limit of \eqref{y} and 
the corresponding controlled system, using a suboptimal control, has also been established. Due to the one-dimensional 
phase space the extension to the hydrodynamic equation is not relevant. However, moments of the kinetic distribution as
well as equilibrium conditions have been studied. 
\par
A further extension to the models has been the addition of stochasticity to the dynamics. In   \cite{Albi:2017aa} white noise, $dW$, is added to the dynamics leading to a problem of the form:
\begin{align*}
&q^* = \mbox{argmin}_{v \in \R} \frac{1}2 \int_0^T \mathbb{E}\left( \left( \frac{1}N \sum\limits_{j=1}^N g(y_i)  \right)^2 + \frac{\beta}2 v^2 \right) dt, \;  \\
& \mbox{ subject  to }  \; dy_i  = \Big(\frac{1}N \sum\limits_{j \not = i }^N P(y_i,y_j) (y_j-y_i)  + v \Big) \;dt + dW_i, \; y_i(0)=y_{i,0},
\end{align*}
where $P$ is a given interaction kernel and $g$ a given cost functional. Also,  in   \cite{Albi:2015ab} a stochastic parameter,
$W$, is included in the collision operator, $P$, to account for modeling errors. Thus the resulting problem 
reads
\begin{align*}
&u = \mbox{argmin}_{v \in \R} \frac{1}2 \int_0^T \mathbb{E}\left( \left( \frac{1}N \sum\limits_{j=1}^N g(y_i)  \right)^2 + \frac{\beta}2 v^2 \right) dt  \mbox{ subject  to }  \eqref{y}.
\end{align*}
Using instantaneous control and polynomial chaos expansion the problem is reduced to a closed loop problem for
suboptimal control $q.$
\par
Additional structural requirements on the control could be considered. For example, some application 
might require sparse control. An existing approach modifies the cost functional 
to treat this case \cite{Fornasier:2014aa}:
\begin{align*}
q^* = \mbox{argmin}_{v \in \R} \frac{1}2 \int_0^T  \frac{1}N \sum\limits_{j=1}^N y_i^2   dt + \beta \| v \|_{L},
\end{align*}
for $L=L_p(0,T)$ where $0 < p \leq 1$ is studied.  This term promotes sparsity of the control in time. 
Another type of sparsity is introduced in \cite{Albi:2016aa} where sparsity is introduced in
the number of controllable particles leading to a problem of the following type:
\begin{align*}
 y_i'  = \frac{1}N \sum\limits_{j \not = i }^N P(y_i,y_j,W) (y_j-y_i)  +  b_i v, \; y_i(0)=y_{i,0},
\end{align*}
where $b_i \in \{ 0, 1\}$, fixed, with $\| b\|_{\ell^1}$ sufficiently small. 

\subsection{Controlled Particle Systems in High--Dimensional State Space}

Typical crowd dynamic models \cite{Bellomo:2013aa,Cristiani:2014aa,Vicsek:2012aa} have a state space that at least contains the velocity $v_i$ 
and the position of the particle $x_i.$ However, the field of control for such models and its associated mean field equations is largely
unexplored. So far, results on instantaneous control approaches applied to second--order alignment  models are available in \cite{Albi:2014ae}. The basic model for $i=1,\dots,N$ particles is given by the following dynamics:
\begin{align}
\frac{d}{dt} x_i(t) &= v_i(t), \\
\frac{d}{dt} v_i(t) &= \frac1N \sum\limits_{j=1}^N P(x_i(t),x_j(t))(v_j(t)-v_i(t)) +q^*(t)
\end{align}
subject to initial conditions and for a given function $P.$ The instantaneous control approach has been used to 
find a suboptimal explicit closed loop control for $q^*(t)$ by minimizing 
\begin{align}\label{eq:Objqstqr}
q^*(t) = \mbox{argmin}_{q \in \R} \frac{1}2 \int_{t}^{t+\Delta t}  \frac{1}N \sum\limits_{j=1}^N (v_i(s)-v_d(s))^2  + \frac{\beta}2  q^2(s) ds ,
\end{align}
on a small time horizon $\Delta t>0.$ The parameter $\beta>0$ is again the regularization parameter and $v_d$ is a desired velocity,
piecewise constant on the receding horizon $(t,t+\Delta t).$ Depending on the kernel, $P$, different crowd dynamics \cite{} 
can be described. To illustrate this point, we consider for parameters $K> 0, \gamma>0, \delta\geq 0$
\begin{align} \label{P CS}
P(x_i,x_j) = \frac{ K }{ (\gamma^2 + \| x_i - x_j\|^2)^\delta }.
\end{align}
This dynamics are known as the Cucker--Smale model. We refer to \cite{Carrillo:2010ac,Cucker:2007aa} for more details on properties as well as motivation. Another example of $P$ is introduced in \cite{Motsch:2014aa} given by 
\begin{align} 
P(x_i,x_j) = \frac{H(|x_i-x_j|)}{\frac{1}N \sum_{j = 1}^N H(|x_i-x_j|)}
\end{align}
where $H$  could be given by \eqref{P CS}, i.e., $H=H(r)= \frac{K}{ (\gamma^2 + r^2)^\delta }.$ An extension of the instantaneous control 
approach to this system is straightforward.
\par
A further class of crowd models where, to the best of our knowledge no control results are available for now, are models of the type
\cite{DOrsogna:2006aa}
\begin{align}
\frac{d}{dt} x_i(t) &= v_i(t), \\
\frac{d}{dt} v_i(t) &=-  \frac1N \sum\limits_{j=1}^N P(x_i(t)-x_j(t)) + (\beta - \gamma |v_i|^2) v_i
\end{align}
where $P(x,y)$ is given as a gradient of a potential, $U$, modeling  interactions,  $P(x_i,x_j) = \nabla_{x_i} U(|x_j-x_i|)$ and 
$\beta,\gamma \geq 0$ are parameters modeling the self propulsion of particles.
\par 
Even more detailed models which have been recently introduced such as in \cite{Carrillo:2010ac} introduce, on the particle level, the concept of a visual cone, restricting the 
possible interaction of particle $i$ with particles $j\not = i$ by geometric conditions. So far, no control results for those models are known.

In applications in finance, there only exist a few results on control of wealth or financial market models. For wealth models of the type \eqref{y} subject to \eqref{eq:Objqstqr}, an instantaneous control approach has been introduced in \cite{during2018kinetic}.
The objective function intends to reduce the variance of wealth among agents and thus the objective function in \eqref{eq:Objqstqr} needs to be replaced by the empirical variance. 
Furthermore, the control of binary wealth interactions leading to a Boltzmann type description has been investigated in \cite{during2018kinetic}.

In a game theoretic setting wealth models have been studied as well. A second order model of the type
\begin{align}\label{wealth}
& \dot{x}_i = V(x_i,y_j),\\
& dy_i = q_i^*\ dt + y_i\ dW_i,
\end{align}
has been discussed in \cite{degond2014evolution}. Here $V(\cdot, \cdot)$ models the speed of change in the economic configuration $x_i\in\R$. The control $q_i^*$ has been computed as best reply to a cost functional of the type
\begin{align}\label{objW}
\displaystyle{q_i^*(t) = \mbox{argmin}_{q \in \R} \frac{1}2  \int_{t}^{t+\Delta t} \mathbb{E}\left(  \frac{1}N \sum\limits_{j=1}^N P(x_i(s), x_j(s))\  \Phi(y_j(s)-y_i(s))  + \frac{\beta}2  q^2(s)\right) ds},
\end{align}
where $\Phi(\cdot)$ models the trading interaction of agents with different wealth levels $x_i>0$. 
A portfolio model of similar structure to \eqref{wealth} -\eqref{objW} has been discussed in \cite{trimborn2019portfolio}, where additionally the dynamics are coupled to a  stock price equation modeled by a stochastic differential equation. \\ \\
An example of a wealth model which has been studied in the game theoretic setting without any sub-optimal strategies is presented in \cite{gueant2011mean}. 
The structure of the model is as follows 
\begin{align} \label{MFSys}
& dx_i = q_i^*\ dt + dW_i,\\
&q_i^*(t) = \mbox{argmin}_{q \in \R} \frac{1}2  \int_{t}^{\infty} \mathbb{E}\left( F(x_1,...,x_N, q_i)\ e^{-r(s-t)} \right)   ds.
\end{align}
In the mean field limit the system \eqref{MFSys} reduces to a Hamilton-Jacobi-Bellmann equation. 
These models are known as mean field games and have been extensively studied in the past decade, see \cite{bensoussan2013mean, carmona2018probabilistic}.

\section{Summary}

In summary chapter presents the some developments in controlled crowd dynamics. A discussion of the multiscale control of particle systems which culminates in the mean field of the controlled particle system. In addition, the ensuing control system for the hydrodynamic model are derived from the particle system. To demonstrate the practical application of the approaches, numerical tests are undertaken on the linear models. The behaviour predicted in the theoretical discussions is clear  in the numerical results. Further a discussion of work that has been done on nonlinear models is also discussed in the last section of the chapter. This includes applications in opinion formation, stochastic control, crowd models as well as wealth models among others.

\begin{acknowledgement}
We acknowledge the support by the National Research Foundation of South Africa (Grant number: 93099 and 93476), DFG HE5386/14,15, 18, BMBF ENets 05M18PAA,  Cluster of Excellence Internet of Production (ID 390621612) and NSF RNMS grant No. 1107291 (KI-Net).
\end{acknowledgement}
%
 \bibliographystyle{siamplain}

 \bibliography{CompleteBibTex}
 
\end{document}